\documentclass[11pt,reqno]{article}
\usepackage[utf8]{inputenc}
\usepackage{color}
\usepackage{mathrsfs}
\usepackage{pifont}
\usepackage{amsmath}
\usepackage{amsthm}
\usepackage{txfonts}
 \usepackage{geometry}
\usepackage{latexsym}
\usepackage{amssymb}
\usepackage{graphicx}%
\usepackage{ifthen}
\usepackage{fancyhdr}
\usepackage{afterpage}
\usepackage{xcolor}
\usepackage{ifthen}

\newlength{\defbaselineskip}
\newcommand{\setlinespacing}[2]%
           {\setlength{\baselineskip}{#1 \defbaselineskip}}
\newcommand{\doublespacing}{\setlength{\baselineskip}%
                           {1.5 \defbaselineskip}}

\newtheorem{lemma}{Lemma}
\newtheorem{theorem}{Theorem}

\newtheorem{cor}{Corollary}

\theoremstyle{definition}

\begin{document}
\doublespacing
\begin{center}
{\Large Sequential Detection of Common Change in High-dimensional Data Stream with Cross-Sectional Dependence} \\
{\it Yanhong Wu{\footnote{Address correspondence to Yanhong Wu, Department of Mathematics, California State University Stanislaus, Turlock, CA 95382, USA; E-mail: ywu1@csustan.edu}} and Wei Biao Wu } \\
Department of Mathematics, California State University Stanislaus,\\
Turlock, California, USA \\
Department of Statistics, University of Chicago, Chicago, Illinois, USA \\
\end{center}

\begin{center}
ABSTRACT
\end{center}
\begin{quote}
{\em After obtaining an accurate approximation for $ARL_0$, we first consider the optimal design of weight parameter for a multivariate EWMA chart that minimizes the stationary average delay detection time (SADDT). Comparisons with moving average (MA), CUSUM, generalized likelihood ratio test (GLRT), and Shiryayev-Roberts (S-R) charts after obtaining their $ARL_0$ and SADDT's are conducted numerically. To detect the change with sparse signals, hard-threshold and soft-threshold EWMA charts are proposed. Comparisons with other charts including adaptive techniques show that the EWMA procedure should be recommended for its robust performance and easy design. }\\
\end{quote}

\noindent{\bf Keywords:} Stationary average delay detection time; Change point; EWMA procedure; MA procedure; Generalized LRT.\\

\noindent{\bf Mathematics Subject Classification:} Primary 62L10; Secondary 62N15

\section{Introduction}
Recent researches on sequential change point detection problems have been focused on developing detecting procedures for multi-stream data or panel data. The first type of models assumes that a change occurs in all or part of the panels, called common change, as considered in Tartakovsky and Veeravalli (2008), Mei (2010), Xie and Siegmund (2013) and Chan (2016). Wu (2019, 2020) and Ning and Wu (2021) proposed a procedure where the sum of Shiryayev-Roberts processes is used as the detecting process and the individual CUSUM processes are used to isolate the changed panels with FDR control after detection. Wu and Wu (2022) also discussed the common change detection in multi-Poisson process case. The second type of models assumes that each panel has it own change point and individual detecting procedures are applied to each panel and terminated when a change is detected.  Chen et. al. (2020) and Chen and Li (2019) considered the detection and isolation of changed panels by treating the problem as multiple sequential hypothesis tests problems and then applying FDR control. A recent review is given by Bardwell, et al. (2020). 
However, most results are given for independent data streams. 

In on-line quality control, a large amount of recent literature are dedicated to multi-variate case  where the  dimension is not very large and the observations may be cross-sectional or serial dependent. Many multivariate CUSUM and EWMA based procedures have been developed for detecting the mean and (or) covariance matrix change. For example, Crosier (1988) and Ngai and Zhang (2001) considered the generalization of CUSUM chart and Lowry, et al. (1992) considered a multivariate EWMA chart for detection of mean.
On the other hand, structural change in time series have been considered from econometric point of view where serial dependence represented as random effects is considered, see Bai (2010), Kim (2014), and Li, et al. (2016),

The basic model can be described as follows. For $i=1,2,...,N$, assume the observations observed  in $i^{th}$ panel are
\[
X_{it} = \mu_{it} + Z_{it},
\]
for $t=1,2,...$, where $\mu_{it}$ is the mean and $Z_{it}$ are the errors that are not only cross-sectionally dependent, but may also be serial-dependent and panel dependent. 
For example, under the simple intra-class model with mean change, $Z_{it} = \frac{1}{\sqrt{N}} a_t +e_{it} $, $\mu_{it} =\mu_{i0}$ for $t \leq \nu $ and $\mu_i $ for $t > \nu$. $\nu$ is the common change point. Here $a_t$ have means 0 and variance $\sigma_a^2$, and  $e_{it}$ have mean 0 and variance $\sigma_e^2$ and all are independent across for both $i$ and $t$. 
More generally, we may generalize $Z_{it} =\gamma_i a_t +e_{it}$ where $\gamma = (\gamma_1, ..., \gamma_p)^T$ represents the direction of correlation such that $\sum \gamma_i^2 =1$. 

In this communication, our main interest is to present some theoretical results for $ARL_0$ and stationary average delay detection time (SADDT) in the multivariate case, where the SADDT is defined as in Srivastava and Wu (1993) and Knoth (2021). In Section 2, we first give the definition of several charts based on EWMA, MA, CUSUM, and Generalized Likelihood Ratio Test that are direction-invariant. 
 In Section 3, we consider the optimal design of MEWMA procedure. Accurate approximation for $ARL_0$ is given for the design of the chart. Optimal value for the weight parameter when the signal strength is known is obtained that minimizes the asymptotic stationary average delay detection time (SADDT). The simple intra-class dependence model and the model with a random loading factor are used as illustrations. Comparisons with other multivariate CUSUM charts are conducted in Section 4 where the theoretical results for $ARL_0$ and $SADDT$ for other charts are also presented. In Section 5, we consider to use the hard or soft threshold technique to reduce the SADDT when signal only appear in a proportion of the channels or is weak in some channels. The use of adaptive estimation when the post-change mean is unknown is also discussed. Our main findings show that the threshold method in EWMA chart works better when the proportion of channels with signal is small. 

\section{Direction-Invariant Multivariate Charts}

We start with the following multivariate change point model with possible cross-dependence. 
Let \[
X_{it} = \mu_i +Z_{it} ,
\]
for $i=1,...,N$ and $t=1,2,...$, where $\mu_i=\mu_i I_{[t > \nu]}$ and $Z_{t} =(Z_{1t},...,Z_{Nt})^T$ are i.i.d. following $N(0, \Sigma)$. $\nu$ is called the common change point. 
By denoting $X_t = (X_{1t}, ..., X_{Nt})^T$, $\mu =(\mu_1,..., \mu_N)^T$, and $||\mu|| =(\mu^T\Sigma^{-1} \mu)^{1/2} $ as the signal strength, we can write $\mu =||\mu|| \Sigma^{-1/2}\gamma $ where $
\gamma = \Sigma^{-1/2} \mu /||\mu|| $ is the signal direction. We see that $X_t$ follows multivariate normal $N(||\mu||\Sigma^{1/2} \gamma  I_{[t>\nu]} , \Sigma )$. For observation $X_t$, given a reference value $||\delta||$ for the signal strength $||\mu||$, the log-likelihood ratio for $t>\nu$ vs $t\leq \nu$ is equal to
\begin{equation}
||\delta||\gamma^T \Sigma^{-1/2} X_t -\frac{||\delta||^2}{2} = ||\delta || (\gamma^T \Sigma^{-1/2} X_t -\frac{||\delta||}{2}) .
\end{equation}
The following several charts have been used to develop a chart that is direction invariant.

\noindent{\bf (1) Multivariate EWMA chart}

The first approach is to use EWMA of the log-likelihoods as discussed in Lowry, et al.(1992) and Wu and Wu (2022). By defining
\[
Y_{t}= (1-\beta)Y_{(t-1)} +\beta X_{t} 
\]
for the same weight parameter $\beta$ with $Y_0=0$, the EWMA of the log-likelihood (2.1) is equivalent to
\[
Y_{t}(\gamma) = \gamma^T\Sigma^{-1/2} Y_t .
\]
If we want the detection procedure is directional invariant, we should use $
\max_{\gamma} Y_t(\gamma)$ as the detection process. This gives $\gamma=\Sigma^{-1/2}Y_t/(Y_t^T\Sigma^{-1}Y_t)^{1/2}$ and $
\max_{\gamma} Y_t^{\gamma} =(Y_t^T\Sigma^{-1}Y_t)^{1/2}$. 

\noindent {\bf \underline{MEWMA Procedure}}: {\em An alarm is made at
\begin{equation}
\tau_{MEW}  =\inf\{ t>0: ~  Y_{t}^T\Sigma^{-1}Y_t >b^2 (\frac{\beta}{2-\beta}) \},
\end{equation}
where the factor $\beta/(2-\beta)$ is the limiting variance for each component. }

\noindent{\bf Remark:} Under the simple independence model without cross-dependence, i,e, $Z_t $ follows $N(0, \sigma_e^2 I)$, the regular MEWMA procedure makes an alarm at
\begin{equation}
\tau_{MEW0} =\inf\{ t>0: ~  \frac{1}{\sigma_e^2} Y_{t}^TY_t >b^2 (\frac{\beta}{2-\beta}) \}.
\end{equation}

\noindent{\bf (2) Multivariate MA chart}

The second approach is to take the moving average of the log-likelihoods with window size, say $w$. By denoting $\bar{X}_{t;w} = \frac{1}{w} \sum_{t-w+1}^t X_j $, the moving average of the log-likelihoods is equivalent to 
\[
\gamma^T\Sigma^{-1/2}\bar{X}_{t;w} .
\]
It is obvious that $\gamma =\Sigma^{-1/2}\bar{X}_{t;w}/ (\bar{X}_{t;w}^T\Sigma^{-1} \bar{X}_{t;w})^{1/2} $ maximizes the above quantity. 

\noindent{\bf \underline{MMA Procedure}}: {\em Make an alarm at 
\[
\tau_{MMA} =\inf\{t>0:  (\bar{X}_{t;w}^T\Sigma^{-1} \bar{X}_{t;w})^{1/2}>h \}. 
\]
}
\noindent{\bf (3) Multivariate CUSUM chart}

There have been considerable discussions on multivariate CUSUM procedure starting with Woodall and Ncube (1985) where individual two-sided CUSUM processes are run simultaneously, and Crosier (1988) where a shrinkage technique for the cumulative signal strength is used.  To develop a directional invariant procedure, we use the following similar argument as for the MEWMA (Ngai and Zhang (2001)).

Note that the log-likelihood ratio for testing $H_0: \nu =\infty$ vs $H_1: 0 \leq \nu < t$ for observations $X_1,.., X_t$ with post-change mean $\mu$ can be written as
\begin{equation}
||\delta||\gamma  \Sigma^{-1/2}(X_{\nu+1}+...+X_t)-\frac{1}{2} (t-\nu)||\delta||^2  .
\end{equation}
The $\gamma$ that maximizes the above log-likelihood ratio is equal to
\[
\Sigma^{-1/2}(X_{\nu+1}+...+X_t) / ((X_{\nu+1}+...+X_t)^T\Sigma^{-1}(X_{\nu+1}+...+X_t))^{1/2}.
\]
Thus, the maximum log-likelihood ratio is equivalent to
\[
((X_{\nu+1}+...+X_t)^T\Sigma^{-1}(X_{\nu+1}+...+X_t))^{1/2} -\frac{k}{2} (t-\nu) .
\]
To reduce the computational complexity, we use the the following window restricted Multivariate CUSUM procedure:

\noindent{\bf \underline{MCUSUM procedure:}} {\em For a large $W$, make an alarm at 
\[
\tau_{MCU} =\inf\{ t>0: \max_{0<w<W} w ((\bar{X}_{t;w}^T\Sigma^{-1} \bar{X}_{t;w})^{1/2}-\frac{k}{2}) > d\}
\]
}

\noindent{\bf Remark:} To mimic the recursive form of one-dimensional process, the following CUSUM procedure was proposed by Pignatiello and Runger (1990), where the common change point estimation is updated at each time $t$ and the cumulative sums after change point estimation are used in the recursive form: 

 {\em Let $\hat{\nu}_0 =0$ and for $t\geq 1$
\begin{equation}
MC1_t =max(0, ((\sum_{\hat{\nu}_t+1}^t X_i )^T \Sigma^{-1} (\sum_{\hat{\nu}_t+1}^t X_i ))^{1/2} - \frac{k_1}{2} (t- \hat{\nu}_t) )
\end{equation}
for a reference value $k_1 $ for the signal strength and reset $
\hat{\nu}_{t+1} =t,$
if $MC1_t =0$; otherwise $\hat{\nu}_{t+1}= \hat{\nu}_t$. An alarm will be made at
\[
\hat{\tau}_1 =\inf\{ t>0: MC1_t > h_1 \}.
\]
After detection, the change point can be estimated as $\hat{\nu}_{\hat{\tau}_1}$. }

To reduce the average delay detection time when the true signal strength is different from the reference value, we may further use the adaptive estimator for the strength of signal. A natural candidate is the EWMA estimator $(Y_t^T\Sigma^{-1} Y_t)^{1/2}$, that can replace $k_1$. Formal discussion of adaptive CUSUM and S-R procedures will be given in Section 5. 

\noindent{\bf (4) Generalized Likelihood Ratio Test (GLRT) based chart}

If we further assume $||\delta||$ is unknown in (2.4), we see the MLE for $||\delta||$ is given by 
$((X_{\nu+1}+...+X_t)^T\Sigma^{-1}(X_{\nu+1}+...+X_t))^{1/2}/(n-\nu) $. Again to control the computational complexity, we can use the following generalized LR based chart. 

\noindent{\bf \underline{GLRT based chart}}:  {\it For a large $W$, make an alarm at
\[
\tau_{GLR} =\inf\{t>w: \max_{0\leq w <W} w (\bar{X}_{t;w}^T\Sigma^{-1} \bar{X}_{t;w}) > b^2 \}.
\]
}

\noindent{\bf Remark.} The window restriction can make the design for the CUSUM chart very complicated. So we shall only consider the window restricted CUSUM chart in the comparison.  The optimality of the Shiryayev-Roberts(S-R) procedure in the one-dimensional case was extended to multivariate case in Siegmund and Yakir (2008) where the window restricted  mixture CUSUM and S-R procedures are studied in terms of $ARL_0$ and maximum of average delay detection time. 

\section{Optimal Design of MEWMA chart} 

In this section, we shall focus on the MEWMA chart. 
For the design of the control limit $b$ and weight $\beta$, the following theorem gives a corrected approximation for $ARL_0 =E_{\infty}[\tau]$: 

\begin{theorem}
As $b\rightarrow \infty $ and $\beta \rightarrow 0$, $\tau_{MEW}$ under $E_{\infty}[.]$ is asymptotically exponential with mean
\begin{eqnarray}
ARL_0 &\approx &  \frac{1}{2\beta} \int_0^{{b^*}^2/2} x^{-N/2}e^x \int_0^x z^{N/2-1}e^{-z} dzdx \nonumber \\
&\approx &  \frac{1}{-2\ln(1-\beta)} \int_0^{{b^*}^2/2} x^{-N/2}e^x \int_0^x z^{N/2-1}e^{-z} dzdx \\
& \approx & \frac{1}{-2\ln(1-\beta)} \Gamma(\frac{N}{2}) (\frac{{b^*}^2}{2})^{-N/2} e^{{b^*}^2/2 }, 
\end{eqnarray}
where $b^*$ is the continuous correction of $b$ given by $b^* =b+\rho_+ \beta /\sqrt{\beta/(2-\beta)}$ with $\rho_+ \approx 0.5826$ (Siegmund (1985, pp. 50)) and $-ln(1-\beta) $ gives better numerical approximations than $\beta$.
\end{theorem}

\noindent{\it Proof.} Standard weak convergence theory shows that as $\beta \rightarrow 0$ and $L \rightarrow \infty$ such that $\beta L\rightarrow \infty$, 
\[
{Y}_j^T \Sigma^{-1} {Y}_j/\sqrt{\beta/(2-\beta)} \Longrightarrow \{e^{-u} ||\textbf{W}(e^{2u})||, 0 \le u \le L \beta \},
\]
where ${\textbf{W}(t)}$ is a $N$-dimensional Brownian motion. Following the proof of Theorem 3  in Wu and Wu (2022), $\tau_{MEW}$ can be treated as the boundary crossing time for a CIR process with drift $\mu(x) =N-2x$ and diffusion $\sigma^2(x) = 4x$. Following the standard notation of diffusion process, since
\[
s(x)= exp\left(-\int^x \frac{2 \mu(u)}{\sigma^2(u)} du \right)=exp\left(-\int^x (\frac{N}{2u} -1) du \right) = x^{-N/2} e^x,
\]
so the scale function $S(x)$ is given by
\[
S(x)=\int_0^x s(u)du = \int_0^x u^{-N/2} e^u du .
\]
On the other hand, the speed function $M(x)=\int_0^x m(u)du$ where the speed density $m(x)$ is given by
\[
m(x) =\frac{2}{\sigma^2(x) s(x) } =\frac{1}{2} x^{N/2-1} e^{-x}. 
\]
The result is obtained from the formula given in Karlin and Taylor (1981, Pg 197 (4.12)) after a boundary correction technique by changing $b$ by $b^* =b+\rho_+ \beta /\sqrt{\beta/(2-\beta)}$ as considered in Wu and Wu (2022). \qed

The simulation results given in Tables 1 and 2 in Example 1 show that the approximation is very accurate.  The advantage of the MEWMA chart is that we do not need to assign a reference value for the change signal strength while all possible directions of change are considered. However, if a reference value for the signal strength $(\delta^T \Sigma^{-1} \delta)^{1/2} $ is given, we can consider the optimal design of the weight parameter $\beta$ that minimizes the stationary average delay detection time as defined in Srivastava and Wu (1993) and Knoth (2021) . The following theorem gives the optimal design for $\beta$ and $b$. 

\begin{theorem}
Assume $||\delta|| $ is the reference value for the signal strength with post-change mean ${\bf \delta} =(\delta_1,..., \delta_N)^T$ such that $0<{\bf \delta}^T\Sigma^{-1}{\bf \delta}<\infty$. Then as $ARL_0 \rightarrow \infty$, the optimal value of $b$ and $\beta$ are given by 
\[
b^* =\sqrt{2\ln (ARL_0)} (1+o(1)) \]
\[
\beta^* = \frac{2{\bf\delta}^T \Sigma^{-1}{\bf \delta}}{{b^*}^2} k^* (1+o(1)) =k^*  \frac{{\bf\delta}^T \Sigma^{-1}{\bf \delta}}{\ln (ARL_0)} (1+o(1)) ,
\]
where $k^*=0.5117$. Under the optimal design, the SADDT is given by
\begin{equation}
SADDT = \frac{ {b^*}^2 c^*} {2 {{\bf \delta}}^T\Sigma^{-1} {\bf\delta} } (1+o(1))
=\frac{c*}{{\bf \delta}^T\Sigma^{-1} {\bf\delta} } \ln(ARL_0 ) (1+o(1)) ,
\end{equation}
where $c^* \approx 2.4554$
\end{theorem}

\noindent{\em Proof.} The results in the one dimensional continuous time Brownian motion case is given in Srivastava and Wu (1993). The proof will be completed after two lemmas.

\begin{lemma}
As $ARL_0 \rightarrow \infty$, $\nu \rightarrow \infty$, the optimal $\beta^* \rightarrow 0$. 
\end{lemma}
\noindent{\em Proof.} By contradiction, if $\beta^* $ does not go to zero, then from (2.2), $b\rightarrow \infty$. In the stationary  state, we can assume $Y_{\nu}$ follows $N(0,\beta /(2-\beta) \Sigma) $. By denoting $\tilde{Y}_t = Y_{\nu+t}$, we can write
\begin{eqnarray}
Y_{\nu+t}^T \Sigma^{-1} Y_{\nu+t} & =& (\tilde{Y}_t + (1-(1-\beta)^t) {\bf \delta})^T \Sigma^{-1} (\tilde{Y}_t + (1-(1-\beta)^t) {\bf \delta}) \nonumber \\
& = & \tilde{Y}_t^T \Sigma^{-1} \tilde{Y}_t +2 (1-(1-\beta)^t ){\bf \delta}^T \Sigma^{-1} \tilde{Y}_t + {\bf \delta}^T\Sigma^{-1}{\bf\delta}(1-(1-\beta)^t )^2  . 
\end{eqnarray}
That means, $Y_{\nu+t}^T \Sigma^{-1} Y_{\nu+t}$ behaves like a non-central stationary chi-square process with center $\delta$. From Theorem 1, we see that $\ln(ARL_0) = b^2/2 (1+o(1)$, This implies that $\ln (SADT) = O(\ln (ARL_0 ) ) $.  \qed 

\begin{lemma}
As $ARL_0 \rightarrow \infty$, $\beta \rightarrow 0$ such that
\begin{equation}
\frac{b^2 \beta }{2{\bf \delta}^T\Sigma^{-1}{\bf \delta}} \rightarrow k <1,
\end{equation}
then
\[
SADDT = -\frac{1}{\beta} \ln (1-\sqrt{k}) (1+o(1)). 
\]
\end{lemma}

\noindent{\em Proof.} From (3.9), as $\beta \rightarrow 0$, we see
\[
\tilde{Y}_t^T \Sigma^{-1}\tilde{Y}_t =_d \frac{\beta}{2-\beta} \chi_N^2,
\]
and
${\bf \delta}^T \Sigma^{-1}\tilde{Y}_t$ follows $\sqrt{\frac{\beta}{2-\beta}} N(0, {\bf \delta}^T \Sigma^{-1}{\bf \delta}) $. So by ignoring the overshoot, at the boundary crossing $t=\tau-\nu  =\tilde{\tau}$, we have
\[
(1-e^{-\beta t})^2 =k -2(1-e^{-\beta t}) \frac{ {\bf \delta}^T \Sigma^{-1} \tilde{Y}_t }{ {\bf \delta}^T\Sigma^{-1}{\bf\delta}} - \frac{ \tilde{Y}_t^T \Sigma^{-1} \tilde{Y}_t}{{\bf \delta}^T\Sigma^{-1}{\bf\delta} }. 
\]
This gives
\[
1-e^{-\beta t} = \sqrt{k} (1- \frac{2}{k} (1-e^{-\beta t}) \frac{ {\bf \delta}^T \Sigma^{-1} \tilde{Y}_t }{ {\bf \delta}^T\Sigma^{-1}{\bf\delta}} - \frac{ \tilde{Y}_t^T \Sigma^{-1} \tilde{Y}_t}{{\bf \delta}^T\Sigma^{-1}{\bf\delta} })^{1/2}
\]
\[
=\sqrt{k} (1- \frac{1}{k} (1-e^{-\beta t}) \frac{ {\bf \delta}^T \Sigma^{-1} \tilde{Y}_t }{ {\bf \delta}^T\Sigma^{-1}{\bf\delta}}+O_p(\beta) )
\]
\[
=\sqrt{k} - \frac{ {\bf \delta}^T \Sigma^{-1} \tilde{Y}_t }{ {\bf \delta}^T\Sigma^{-1}{\bf\delta}}+O_p(\beta).
\]
Thus, 
\[
-\beta t = \ln(1-\sqrt{k} - \frac{ {\bf \delta}^T \Sigma^{-1} \tilde{Y}_t }{ {\bf \delta}^T\Sigma^{-1}{\bf\delta}} +O_p(\beta)) 
\]
\[
=\ln(1-\sqrt{k})   -\frac{1}{1-\sqrt{k}} \frac{ {\bf \delta}^T \Sigma^{-1} \tilde{Y}_t }{ {\bf \delta}^T\Sigma^{-1}{\bf\delta}}  +O_p(\beta).
\]
So $SADDT =E_0^*(\tau_{MEW}) = -\ln(1-\sqrt{k})/\beta$. We can also obtain the first order variance.  \qed

Since $\beta = 2k \delta^T\Sigma\delta /b^2$, we have
\begin{equation}
SADDT = \frac{b^2}{2{\bf \delta}^T\Sigma^{-1}{\bf \delta}} \frac{\ln (1-\sqrt{k})}{-k } (1+o(1)).
\end{equation}
The proof of Theorem 2.2 is completed by noting that $\frac{\ln (1-\sqrt{k})}{-k }$ reaches its minimum value $c^* =2.4554$ at $k^* =0.5117$. \qed

By checking the proofs for the lemmas, we can derive the SADDT for any signal strength $||\mu|| \neq ||\delta||$ that will be useful when we study the efficiency of the EWMA procedure compared with other charts. The following corollary gives the result for SADDT that only depends on the strength of signal. 

\begin{cor}
Suppose the true post-change mean is ${\bf \mu} =(\mu_1, ..., \mu_N)^T$ with signal strength $\mu^T \Sigma^{-1} \mu$ and $\beta^* $ and $b^*$ are selected as in Theorem 2. Then if
\[
k= k^* \frac{{\bf \delta}^T \Sigma^{-1}{\bf \delta} }{{\bf \mu}^T \Sigma^{-1}{\bf \mu} } <1 ,
\]
the MEWMA procedure is efficient and the SADDT is asymptotically equal to
\[
SADDT = \frac{\ln(ARL_0)}{{\bf \mu}^T\Sigma^{-1}{\bf \mu}} \frac{\ln (1-\sqrt{k})}{-k } (1+o(1)) = \frac{\ln (ARL_0)}{\delta^T\Sigma^{-1}\delta} \frac{\ln (1-\sqrt{k^*} (\frac{\delta^T\Sigma^{-1}\delta}{\mu\Sigma^{-1}\mu})^{1/2})}{-k^*} (1+o(1)).
\]
\end{cor}

\noindent{\bf Remark}: The above results are given for finite dimension $N$. When $N \rightarrow \infty$, the results are still true by assuming $N\ln(N) =o(\ln (\beta ARL_0))$ as $\ln \Gamma (N/2)) \approx (N/2) (\ln (N/2/e) -1)$. We also assume the signal strength ${\bf \mu}^T\Sigma{\bf \mu}< \infty$. If the signal strength goes to zero, further discussion will be necessary.

\noindent{\bf Remark}: We used the term efficiency that means the stationary average delay detection time should be at the order of $\ln (ARL_0)$.

\noindent{\bf Example 1: Intra-class dependence model}

Let
\[
X_{it} = \mu_i +\frac{1}{\sqrt{N}} a_t +e_{it} ,
\]
for $i=1,...,N$ and $t=1,2,...$, where $\mu_i=\delta_i I_{[t>\nu]}$, $a_t$ are i.i.d following $N(0, \sigma_a^2 )$, and $e_{it}$ are i.i.d. following $N(0, \sigma_e^2)$. Also, $a_t$ and $e_{it}$ are also independent. 

By denoting $X_t = (X_{1t}, ..., X_{Nt})^T$, $\mu =(\mu_1,..., \mu_N)^T$, $\theta = \sigma_a^2 /\sigma_e^2$ and $\rho= \theta /(1+\theta)$, we can see $X_t$ follows multivariate normal $N(\mu_0 +\delta I_{[t>\nu]} , \Sigma )$ 
with 
\[
\Sigma =\sigma_e^2 (I + \frac{\theta}{N} J),
\]
where $I_{N\times N}$ is the identity matrix, and $J_{N\times N}$ is the matrix with elements being 1. 

From Sylvester' theorem, we know
\[
\Sigma^{-1} = \frac{1}{\sigma_e^2 } (I-\frac{\theta}{N(1+\theta)} J) =\frac{1}{\sigma_e^2 } (I- \frac{\rho}{N} J),
\]
and 
\[
det(\Sigma) = (\sigma_e^2)^N (1+\theta) = \frac{(\sigma_e^2)^N}{1-\rho} .
\]

For each panel, assume $\mu_{i0} =0$, we calculate the EWMA process by
\[
Y_{it}= (1-\beta)Y_{i(t-1)} +\beta X_{it} 
\]
for the same weight parameter $\beta$ and denote $Y_t = (Y_{1t},..., Y_{Nt})^T$. For the control limit $b^2$, an alarm will be raised at
\begin{eqnarray}
\tau &= & \inf\{t>0: ~ {Y}_t^T \Sigma^{-1} {Y}_t >b^2 (\frac{\beta}{2-\beta})\}  \nonumber\\
&=& \inf\{t>0: ~ \frac{1}{\sigma_e^2} (\sum_{i=1}^N Y_{it}^2 -\frac{\rho}{N} (\sum_{i=1}^N Y_{it})^2 ) > b^2 (\frac{\beta}{2-\beta})\}.
\end{eqnarray}

The following corollary gives the stationary average delay detection time. 

\begin{cor}
Under the intra-class model, for the signal $\mu$, assume $\beta \rightarrow 0 $ and $b \rightarrow \infty$ such that
\begin{equation}
\frac{\beta b^2 \sigma_e^2}{ 2(\mu^T\mu -\frac{\rho}{N} (\mu^T\textbf{1})^2 )} \rightarrow k <1.
\end{equation}
Then the signal $\mu $ is detectable with stationary average delay detection time
\[
SADT = \frac{b^2 \sigma_e^2}{ 2(\mu^T\mu -\frac{\rho}{N} (\mu^T\textbf{1})^2 )} \frac{\ln(1-\sqrt{k})}{-k} (1+o(1)) .
\]
\end{cor}

The introduction of the intra-class model can be used to check how the cross-correlation affects the performance of the regular control chart $\tau_0$ given by (2.3).
For N=10, Table 1 gives the simulated results when $ARL_0 =1000$,$\nu=100$, $\beta =0.01,0.05,0.10$, and $\sigma_a^2= 0, 1, 2, 3, 4, 5$ and $\sigma_e^2=1,0.9,0.8,0.7,0.6,0.5$ correspondingly to make $\sigma_e^2 +\frac{\sigma_a^2}{10} =1$. The control limit $b^*$ is designed from the approximation (2.2) as 4.73 for $\beta =0.01$ and is corrected to 4.64 by using $b =b^*-\rho_+ \beta/ (\beta/(2-\beta))^{1/2}$. Similarly, the values of $b$ for $\beta =0.05$ and 0.10 are $5.14$ and 5.276, respectively. 

With 50,000 replications, reported are the false alarm probability (FAR) and average delay detection time given the change is detected, along with the simulated $ARL_0$. We see that the intra-class coefficient dramatically decreases the $ARL_0$ and thus increases the false alarm probability. However, the average delay detection time keeps quite stable as the mean only changes marginally.

\begin{table}
\begin{center}
Table 1. Simulated $ARL_0$, FAR and SADDT for $N=10$, $\nu=100$ and $ARL_0=1000$  with $\tau_0$\\

{\small
\begin{tabular}{|c |c |c |c | c| c|c |c|}  \hline
$\mu$ & $\beta $ &$\theta=0$ & 1/0.9 & 2/0.8  & 3/0.7 & 4/0.6  &  5/0.5  \\ \hline
0.0 & 0.01 &$ARL_0=1023.0$ & 874.2 & 682.1 & 542.9 & 448.2 & 389.7 \\
    &0.05&  1055.1 & 776.2 & 442.4 & 308.5 & 245.2 & 215.6 \\ 
    &0.1 & 1028.7 & 674.8 & 383.8 & 248.5 & 203.7 & 85.9 \\\hline
0.2 &0.01 & 118.2 & 118.1 & 122.2 & 120.4 & 121.7 & 117.7 \\
    &   & 0.0236 & 0.0292 & 0.0496 & 0.0986 & 0.1392 & 0.1836 \\
    &0.05 & 181.3 & 130.3 & 108.2 & 90.5 & 80.8 & 72.3 \\
     &  & 0.0772 & 0.1078 & 0.2118 & 0.3222 & 0.4466 & 0.537 \\
     &0.1 & 288.2& 165.4 & 111.9 & 85.0 & 65.0& 54.9  \\
     &  & 0.0858 & 0.1342 & 0.2922 & 0.4526 & 0.608 & 0.7124 \\ \hline
0.3  &0.01 &  50.34 & 51.90 & 53.09 & 55.10 & 55.63 & 56.67 \\
     &  & 0.0182 & 0.0296 & 0.0522 & 0.0896 & 0.1334 & 0.1818 \\
    & 0.05 & 43.95 & 41.27 & 40.73 & 39.92 & 38.11 & 37.78 \\
     &  &0.073 & 0.1118 & 0.2188 & 0.331 & 0.4452 & 0.5468 \\ 
     &0.1 & 59.37 & 47.80 & 41.30 & 38.53 & 34.95 & 31.84 \\
     &  & 0.081 & 0.1422 & 0.2852 & 0.4546 & 0.603 & 0.7118\\ \hline
0.5 & 0.01 & 23.05 & 23.76 & 24.77 & 25.37 & 26.17 & 26.74 \\
   &    &  0.0204 & 0.0246 & 0.0528 & 0.1016 & 0.1352 & 0.1806 \\ 
    &0.05 &  14.76 & 15.09 & 15.89 & 16.18 & 16.17 & 16.27 \\
    &   & 0.085& 0.1142 & 0.2118& 0.3314& 0.4542 &0.5448 \\
    &0.1 & 13.59 & 13.92 & 14.01 & 13.89 & 14.62 & 14.04 \\
    &   &0.0928 & 0.1382 & 0.2844& 0.4552 & 0.5952 & 0.7124 \\ \hline
0.6& 0.01 &15.16 &15.60 &16.27 & 16.87 &17.06 &17.66 \\
    &   &0.0182 & 0.0262 & 0.0498 & 0.0958 & 0.1372 & 0.1816 \\
    &0.05 & 8.98 &9.39 & 9.56 & 9.91 & 10.09 & 10.38 \\
    &   & 0.0788 & 0.108 & 0.2152 & 0.327 & 0.458 & 5546 \\
    &0.1 & 7.58 & 7.85 & 7.97 & 8.15 & 8.66 & 8.37 \\
    &   &0.088 & 0.1376 & 0.278 & 0.4502 & 5998 & 0.6986 \\ \hline
0.8 & 0.01 & 11.43 & 11.65 & 12.08 & 12.42 & 12.67 & 13.00 \\
    &   &0.0214 &0.0328 & 0.0524 & 0.0848 & 0.1376 & 0.1864 \\
    &0.05 & 6.54 & 6.77 & 6.91 & 7.27 & 7.23 & 7.47 \\
    &   &0.0818 & 0.113 & 0.221 & 0.3302 & 0.4508 & 0.5468 \\
    &0.1 & 5.29 & 5.45 & 5.67 & 5.82 & 5.87 & 6.14 \\
    &   &0.0774 & 0.1402 & 0.2742 & 0.4674 & 0.6124 & 0.7174 \\ \hline
0.9 & 0.01 & 10.02 & 10.33 & 10.79 & 11.08 & 11.27 & 11.72 \\
    &   &0.0192 & 0.024 &0.0542 & 0.0978 & 0.1388 & 0.1838 \\
    &0.05 & 5.73 & 5.94 & 6.13 & 6.34 & 6.39 & 6.60 \\
    &   &0.0744 & 0.1052 & 0.1956 & 0.3372 & 0.4368& 0.546 \\
    &0.1 & 4.63 & 4.80 & 4.93 & 5.08 & 5.19 & 5.42 \\
    &  &0.0862 & 0.135 & 0.2742 & 0.455 & 0.6078 & 0.7192 \\ \hline
1.0 & 0.01 & 9.06 & 9.41 & 9.63 & 9.87 & 10.26 & 10.45 \\
    &   & 0.0194 &0.0258 & 0.0574 & 0.1002 & 0.1408 & 0.1936 \\
    &0.05 & 5.17 & 5.29 & 5.49 & 5.62 & 5.74 & 5.92 \\
    &   & 0.0764 & 0.1178 & 0.2098 & 0.3328 & 0.4528 & 0.5584 \\
    &0.1 & 4.14 & 4.32 & 4.38 & 4.53 & 4.66 & 4.79 \\
    &   &0.0928 & 0.1422 & 0.2814 & 0.446 & 0.6122 & 0.722 \\ \hline 
\end{tabular} 
}
\end{center}
\end{table}

Table 2 gives the corresponding results by using the control chart $\tau$ given by (3.12) by considering the intra-class correlation coefficient $\rho$ . We see that the false alarm probabilities are roughly the same. However, as the intra-class correlation coefficient increases, the average delay detection time increases. Also, we see the approximation (2.2) gives very satisfactory approximation for $ARL_0$ for all values of $\beta$ with continuous boundary correction.

\begin{table}
\begin{center}
Table 2. Simulated $ARL_0$, FAR and ADT for $N=10$, $\nu=100$ and $ARL_0=1000$ with $\tau_{MEW}$\\

{\small
\begin{tabular}{|c | c |c |c | c| c|c |c|}  \hline
$\mu$ & $\beta $& $\theta=0$ & 1/0.9 &2/0.8  & 3/0.7 & 4/0.6  &  5/0.5  \\ \hline
0.0 & 0.01 &$ARL_0=990.8$ & 1023.4& 1001.1 & 983.3 & 999.3 & 1001.0 \\
    &0.05 &  1035.5 & 1079.5 & 1040.6 & 1051.9 & 1027.6 & 1018.8 \\
    &0.1 & 1049.1 & 1023.4 & 1016.1 & 1008.5 & 997.1 & 1006.7 \\  \hline
0.2 &0.01& ADT= 117.3 & 181.7 & 232.0 & 287.6 & 326.9 & 359.0 \\
     &  & FAR= 0.0188 & 0.019 & 0.0186 & 0.0184 & 0.02 & 0.0188 \\
     & 0.05 & 183.1 & 332.6 & 431.5 & 512.4 & 556.2 & 627.8 \\
     &  &0.0744 & 0.0726 & 0.0792 & 0.0716 & 0.0732& 0.083 \\
     &0.1 & 297.6 & 464.7 & 574.1 & 651.2 & 701.6 & 731.2 \\
     &  & 0.0892 & 0.0872 & 0.0914 & 0.0838 & 0.0838 & 0.0864 \\ \hline
0.3  &0.01 &  49.28 & 74.33 & 93.17 & 111.71 & 128.12 & 146.38 \\
    &   & 0.0172 & 0.0164 & 0.0188 & 0.0168 & 0.0184 & 0.0178 \\
     &0.05 & 43.79 & 83.67 & 127.59 & 167.55 & 214.13 & 246.68 \\
     &  &  0.0754 & 0.0762 & 0.075 & 0.07 & 0.076 & 0.0698 \\
     &0.1& 59.58 & 131.92 & 215.1 & 273.7 & 331.5 & 377.8 \\
     &  &0.0888& 0.0862 & 0.0864 & 0.0916 & 0.0908 & 0.0838 \\ \hline
0.5 & 0.01 & 23.21 & 32.92 & 40.86 & 47.63 & 54.09 & 60.74 \\
   &    & 0.0194 & 0.016 & 0.0208 & 0.0186 & 0.0184 & 0.0204 \\
    &0.05 &  14.86 & 23.53 & 31.78 & 41.03 & 49.45 & 58.31 \\
    &   & 0.0802 & 0.073 & 0.0694 & 0.083 & 0.0758 & 0.0738 \\ 
    &0.1 & 13.62 & 24.86 & 38.15 & 53.04 & 70.17 & 88.58 \\
    &   & 0.081 & 0.086 & 0.0784 & 0.0882 & 0.0844& 0.0936 \\ \hline
0.6&0.01 & 15.05 & 21.13 & 26.12 & 30.14 & 34.20 & 37.35 \\
    &   &0.0182 &0.0154 &0.0212 & 0.0172 & 0.0188 & 0.018 \\
    &0.05 & 8.99 & 13.18 & 17.21 & 21.08 & 24.60 & 28.56 \\
    &  & 0.0786 & 0.08 & 0.075 & 0.071 & 0.0768 & 0.0742 \\
    &0.1 &7.61 & 12.06 & 16.61 & 21.53 & 27.17 & 32.93 \\
    &   &0.0866 &0.0892 & 0.792 & 0.0826 & 0.0938 & 0.0846 \\ \hline
0.8 & 0.01 & 11.21 & 15.61 & 19.10 & 22.21 & 24.88 & 27.41 \\
    &   & 0.0166 & 0.0176 & 0.0176 & 0.019 & 0.0152 & 0.02 \\
    &0.05 & 6.51& 9.32 & 11.91 & 13.93 & 16.33& 18.28 \\
    &   &0.079 & 0.0714 & 0.075 & 0.0766 & 0.0736 & 0.0786 \\
    &0.1   &5.31 & 7.88 & 10.39 & 13.06 & 15.33 & 18.02 \\
    &   & 0.0859 &0.084 & 0.084 & 0.084 & 0.0886 & 0.0844 \\ \hline
0.9 &0.01 & 10.08 & 13.88 & 16.95 & 19.31 & 21.97 & 24.16 \\
    &   & 0.0186 & 0.0164 & 0.0182 & 0.0162 & 0.0166 & 0.0182 \\
    &0.05 & 5.72 & 8.15 & 10.25 & 12.22 & 13.93 & 15.75 \\
    &   & 0.0742 & 0.0728 &0.778 & 0.0744& 0.0778 & 0.081 \\
    &0.1    &4.61 & 6.73 & 8.74 & 10.68 & 12.82 & 14.98 \\ 
    &   &0.0774 & 0.0856 & 0.0794 & 0.0832& 0.0894 & 0.0864 \\ \hline 
1.0 &0.01 &9.11 & 12.47 & 15.23 & 17.62 & 19.78 & 21.62 \\
    &   &0.0166 & 0.0182 & 0.0186 & 0.0174 & 0.0176 & 0.017 \\
    &0.05& 5.18 & 7.25 & 9.01 & 10.73 & 12.29 & 13.54 \\
    &   & 0.078 & 0.0688 & 0.074 & 0.079 & 0.0786 & 0.079 \\
    &0.1  &4.12 & 6.00 & 7.58 & 9.18 & 10.82 & 12.35 \\
    &   &0.0906 & 0.0854 & 0.0822& 0.0884 & 0.0858 & 0.0934 \\ \hline 
\end{tabular} }
\end{center}
\end{table}

\noindent{\bf Example 2: Random factor of loading model}

The intra-class model treats all panels as equally correlated. This may limit its applications as the correlations may depend on the location of panel or the individual conditions. Here, we generalize the model with a random factor of different loading. Let
\[
X_{it} = \delta_i I_{[t>\nu]} +\gamma_i a_t +e_{it} ,
\]
for $i=1,...,N$ and $t=1,2,...$, where $\gamma =(\gamma_1,.., \gamma_N)^T$ is the unit loading that may represent the effect of each individual panel such that $\gamma^T\gamma =1$. $X_t$ follows multivariate normal $N(\mu_0 +\delta I_{[t>\nu]} , \Sigma )$ with
\[
\Sigma =\sigma_e^2 ( I + \theta \gamma \gamma^T). 
\]
From Sylvester's theorem, we know
\[
\Sigma^{-1} = \frac{1}{\sigma_e^2 } (I- \frac{\theta}{1+\theta} \gamma \gamma^T) =\frac{1}{\sigma_e^2 }(I-\rho\gamma \gamma^T)   ,
\]
and 
\[
det(\Sigma) = (\sigma_e^2)^N (1+\theta) =(\sigma_e^2)^N \frac{1}{1-\rho}.
\]

For the control limit $b^2$, an alarm will be raised at
\begin{eqnarray*}
\tau &= & \inf\{t>0: ~ {Y}_t^T \Sigma^{-1} {Y}_t >b^2 (\frac{\beta}{2-\beta})\} \\
&=& \inf\{t>0: ~ \frac{1}{\sigma_e^2} (\sum_{i=1}^N Y_{it}^2 -\rho (\sum_{i=1}^N \gamma_i Y_{it})^2 ) > b^2 (\frac{\beta}{2-\beta})\}.
\end{eqnarray*}

The following corollary gives the first order result for the stationary average delay detection time and its proof is similar to the intra-class model.

\begin{cor}
For the signal $\mu$, assume $\beta \rightarrow 0 $ and $b \rightarrow \infty$ such that
\begin{equation}
\frac{b^2 \beta \sigma_e^2}{ 2 (\mu^T\mu -\rho (\mu^T\gamma)^2 )} \rightarrow k <1. 
\end{equation}
Then the signal $\mu $ is detectable with stationary average delay detection time
\[
SADT = \frac{b^2 \sigma_e^2}{ 2 (\mu^T\mu -\rho (\mu^T\gamma)^2 )}  \frac{\ln(1-\sqrt{k})}{-k} (1+o(1)) .
\]
\end{cor}

\section{Comparison of Multivariate Charts}
In this section, we compare the MEWMA charts with other multivariate charts in terms of average delay detection time. Without loss of generality, we shall assume $\Sigma =I_N$. 

We first give approximations for $ARL_0$ and $SADDT$ for other multivariate charts.

\begin{theorem}
{\it For the MMA chart with window width $w \rightarrow \infty$,
\[
ARL_0 \approx \frac{w \Gamma(N/2)}{2 (h^2w/2)^{N/2}} e^{h^2w/2+\sqrt{2} \rho h} .
\]
For given reference signal strength $\delta^T\Sigma^{-1}\delta$ and $ARL_0$, the optimal $h$ and $w$ are given by
\[
h^* =(\delta^T\Sigma^{-1}\delta)^{1/2} (1+o(1)),~~~~{\rm and } ~~~~ w^*= \frac{2}{\delta^{-1}\Sigma^{-1}\delta} \ln (ARL_0) (1+o(1)),
\]
and under the optimal design
\[
SADDT \approx w^* = \frac{2}{\delta^T\Sigma^{-1}\delta} \ln (ARL_0) (1+o(1)).
\]
}
\end{theorem}

\noindent{\it Proof:} The approximation of $ARL_0$ can be obtained by using the method given in Siegmund, et al. (2010) and Wu and Siegmund (2022). 

As $w\rightarrow \infty$,
\[
w = \frac{2}{h^2} \ln (ARL_0) (1+o(1)).
\]
At the reference signal $\delta$, after the change-point $\nu$, 
\[
E_{\nu}[\bar{X}_{t+\nu-w, t+\nu}^T\Sigma^{-1} \bar{X}_{t+\nu-w, t+\nu}] \approx \delta^T\Sigma^{-1}\delta (\min(t,w))^2/w^2 .
\]
Therefore, in order the signal to be detectable, 
\[
\delta^T\Sigma^{-1}\delta \geq h^2.
\]
Under this condition, 
\[
SADDT = \frac{h w}{(\delta^T\Sigma^{-1}\delta )^{1/2}} (1+o(1)).
= \frac{2}{h (\delta^T\Sigma^{-1}\delta )^{1/2}} \ln (ARL_0) (1+o(1)).
\]
Thus, to minimize SADDT, $h^*=(\delta^T\Sigma^{-1}\delta )^{1/2} $, and 
\[
w^* =\frac{2}{\delta^T\Sigma^{-1}\delta} \ln (ARL_0) (1+o(1)),
\]
and 
\[
SADDT = \frac{2}{\delta^T\Sigma^{-1}\delta} \ln (ARL_0) (1+o(1)).\qed
\]

\begin{cor}
Under the optimal design for MMA procedure as in Theorem 3, if the true signal is $\mu=(\mu_1, ..., \mu_N)^T$ such that $\mu^T\Sigma^{-1}\mu \geq \delta^T \Sigma^{-1}\delta$, then 
\[
SADDT = \frac{2}{\mu^T\Sigma^{-1}\mu} \ln (ARL_0) (1+o(1)). 
\]
When $\mu=(\mu_1, ..., \mu_N)^T$ such that $\mu^T\Sigma^{-1}\mu < \delta^T \Sigma^{-1}\delta$, MMA becomes inefficient. 
\end{cor}

The optimality of MCUSUM procedure is studied in Basseville and Nikiforov (1993) and Nikiforov (1999). Here we summarize the results in the following theorem.

\begin{theorem}
{\it For the window restricted MCUSUM chart with reference signal strength $k =||\delta|| $, as $d \rightarrow \infty$ and $W \rightarrow \infty$ such that $W/d > 2/||\delta|| $, 
\[
ARL_0 = (4||\delta||d)^{-(N-1)/2} \frac{\Gamma(N-1)}{\Gamma((N-1)/2)}\frac{2}{||\delta||^2} e^{||\delta|| (d+2\rho_+)} (1+o(1)),
\]
and when the true signal is $\mu$ such that $||\mu|| >||\delta||/2$, 
\[
SADDT =\frac{1}{||\delta|| (||\mu||-||\delta||/2)} \ln ARL_0 (1+o(1)).
\]
}
\end{theorem}

\noindent{\it Proof:} The approximation for $ARL_0$ can be obtained from Wu and Siegmund (2022). From this, we have
\[
d =\frac{\ln (ARL_0)}{||\delta||} (1+o(1)).
\]
At the boundary crossing time $t >\nu$, for a window size $w$, 
\begin{eqnarray*}
& & \max w((\bar{X}_{t;w}^T\Sigma^{-1}\bar{X}_{t;w})^{1/2} -||\delta||/2) \\
& & =\max  w(((Z_{t,w}/\sqrt{w} +\min(\frac{t-\nu}{w},1)\mu)^T\Sigma^{-1} (Z_{t,w}/\sqrt{w} +\min(\frac{t-\nu}{w},1)\mu))^{1/2}-||\delta||/2)\\
& & = \max w( ((\min(\frac{t-\nu}{w},1))^2 \mu^T\Sigma^{-1}\mu +2 \min(\frac{t-\nu}{w},1) 
\frac{\mu^T\Sigma^{-1}Z_{t;w}}{\sqrt{w}} + \frac{Z_{t;w}^T\Sigma^{-1}Z_{t;w}}{w})^{1/2} -
||\delta||/2) =d. 
\end{eqnarray*}
Thus, as $d\rightarrow \infty$, $w\rightarrow \infty$. So the second and third terms inside the square root are at lower order. Therefore, at the first order, we have
\[
\max w( ||\mu|| \min(\frac{t-\nu}{w},1) - ||\delta||/2) +O_p(1/\sqrt{w})  = d.
\]
Obviously, the maximum value for $w$ is $t-\nu$. This gives
\[
(t-\nu) (||\mu||-||\delta||/2) +O_p(1/\sqrt{w}) = d = \ln (ARL_0)/||\delta|| (1+o(1)),
\]
i.e. $SADDT = d/(||\mu||-||\delta||/2) (1+o(1))$.  \qed

Finally for the GLRT chart, we have the following the results and the proof is given in the appendix. 

\begin{theorem}
{\it For the GLRT chart, as $b \rightarrow \infty$ and $W \rightarrow \infty$ such that $b/\sqrt{W} =O(1)$,
\[
ARL_0 = \frac{\Gamma(N/2)}{2} (b^2/2)^{-N/2} e^{b^2/2} /\int_{b/\sqrt{W}} ^{\infty} \frac{u\upsilon^2(u)}{2} du,
\]
where $\upsilon(u) \approx \exp(-\rho_+ u)$ and for any true signal $\mu$, 
\[
SADDT =\frac{2}{||\mu||^2} \ln ARL_0 (1+o(1)).
\]
}
\end{theorem}

Comparing the results of Theorems 2 to 5, we can see that at the first order, GLRT chart can serve as the benchmark to define the relative efficiency of other charts by taking the ratio of the SADDT of GLRT chart with the reference chart. 

In Table 3, we conduct a simulation comparison between EWMA, MA, CUSUM, and GLRT charts.

For $||\mu|| = 0.25, 0.5, 0.75, 1.0, 1.25, 1.5,1.75, 2.0$, we select $ARL_0 =1000$ and $N=20$. We consider two extreme scenarios: (i) $\mu = ||\mu|| (1, 0, ..., 0)$ and (ii) $\mu =||\mu|| (1/\sqrt{20}, ..., 1/\sqrt{20})$ to make the two have the same signal strength.

(i) For MEWMA chart, we take $\beta =0.05$ and
$b^2 \beta /(2-\beta) =1.07$; 

(ii) for MMA chart, we take $w=20$ and $h^2 =2.1125$; 

We also include window restricted CUSUM  and the GLRT charts with window size $w=20$. That means, we make an alarm at 
\[
\tau_{MCU} =\inf\{ t>0: \max_{1\leq w \leq 20} w ((\bar{X}_{t;w}^T\Sigma^{-1} \bar{X}_{t;w})^{1/2}-\frac{k}{2}) > d\}
\]
and 
\[
\tau_{GLR}=\inf\{t>w: \max_{1\leq w \leq 20} w (\bar{X}_{t;w}^T\Sigma^{-1} \bar{X}_{t;w}) > b^2 \}. 
\]

(iii) For the GLRT chart, $b=7.08$;

(iv) for the CUSUM chart, we select $k=0.5$ and $d \approx 24.15$; 

(v) for the recursive CUSUM from (2.5), we select $k=0.5$ and $d \approx 31$.  

As a benchmark, we also include the Shiryayev-Roberts (S-R) procedure with the same reference value for all channels. It makes an alarm at
\[
\tau_{MSR} = \inf\{ t>0: R_t >B \}
\]
where $R_t = (1+R_{t-1})\exp(\delta \sum_{i=1}^N (X_{it} -\delta/ 2))$. Due to the martingale property, the $ARL_0$ for $\tau_{MSR} $ can be approximated by $B\exp( \delta \sqrt{N} \rho_+)$.

(vi) for S-R chart, we take $\delta =0.5/\sqrt{20}$ and alarm is made at
\[
\tau_{SR} =\inf\{ t>0: R_t = (1+R_{t-1}) \exp (-(0.5/\sqrt{20}) \sum_{i=1}^{20} (X_{it}- (0.5/\sqrt{20})/2) >B \}.
\]
This is equivalent to a one-dimensional S-R procedure with reference value 0.5. From the approximation $ARL_0 \approx B e^{0.5\rho_+}$, we get $B= 747.29$.  

In addition, we also add a scenario for S-R procedure for the harmonic signal with direction proportional to $(1,1/2,..., 1/20)$. All simulations are replicated 10,000 times. The major findings are:

First, all charts are direction-free. The S-R chart which is highly dependent on the direction. Second, the recursive CUSUM chart performs better only when the signal strength is very weak ($||\mu||= 0.25$) as we take the reference value for signal strength as 0.5. However, it becomes less favorable when the signal strength gets larger and also it has quite large false alarm probability.  The EWMA chart performs better for moderate signal strength ($0.5\leq ||\mu || \leq 1.5$). The windowed GLRT chart performs better for large signal strength, as SADDT becomes small. However, its SADDT is much larger for smaller signal strengths.  Third, EWMA chart has the smallest false alarm probability. Overall, the EWMA chart should be recommended due to its simple design, computational convenience and overall stable performance.

\begin{table}
\begin{center}
Table 3.  Comparison of SADDT for $ARL_0 =1000$, $N=20$ and $\nu=100$ \\
\begin{tabular}{|c |c |c |c | c| c|c |}   \hline
$\mu$ & EWMA & MA & GLRT & CUSUM &R-CUSUM & S-R \\ \hline
$ARL_0$ & 1020.5 & 1048.96 & 1010.5& 1032.1 &1024.58 & 1001.95  \\
FAR & 0.0704 & 0.0764 &0.0908& 0.0772&  0.119& 0.0717   \\ \hline
& Scenario (i) &  &    &   &  & \\
0.25 & 382.65 & 578.39 & 705.58 &569.74& 349.488& 473.22   \\
0.5 & 93.65 & 172.78 & 296.37 &174.51 & 106.45&  251.82 \\
0.75 & 41.22 & 56.85 & 99.01 & 56.74& 60.76 &  149.05  \\
1.00 & 25.09 & 27.47 & 37.88 &26.92& 42.26 &  96.63  \\
1.25 & 18.11 & 18.13 & 18.91 &17.25& 33.39&  67.40 \\
1.50 & 14.06 & 14.49 & 12.93 &13.91&  27.43&  50.02  \\
1.75 & 11.56 &  12.48 & 9.76 & 11.62 & 23.11 &  39.77   \\
2.00 & 9.86  &  11.01 & 7.62 & 10.33 & 19.71 & 31.92   \\ \hline 
  & Scenario (ii) &  &    &   &  & \\ 
0.25 &    387.76  &579.15   & 714.15& 607.16 & 345.09&  80.12(161.72) \\
0.5 & 94.49  &170.38  &308.73& 164.07& 107.16&  27.12(54.63)   \\
0.75 & 40.84  & 56.66 & 100.70 & 59.45& 60.63 & 15.70(29.59)  \\ 
1.00 & 25.06 & 27.54  & 38.13 &26.29 & 42.68&  11.10(19.60)  \\ 
1.25 & 18.04  & 18.01   &19.64 & 17.61& 32.61& 8.48(14.58)  \\
1.50 & 14.12  & 14.46  & 13.08 & 13.94 &27.26&  6.92(11.66)    \\
1.75 & 11.55  & 12.45  & 9.45 &11.84 & 22.56&   5.94(9.75)   \\
2.00 & 9.81 &  11.04 & 7.66& 10.21 & 20.13 &   5.14 (8.28)  \\ \hline 
\end{tabular}
\end{center}
\end{table}

\section{Reducing SADDT by Using Threshold Technique and Adaptive Estimation}

\subsection{ By threshold technique}

To reduce the SADDT for detecting multi-dimensional signals, particularly, when signal only appears in a small portion of the dimensions, there are several techniques available. For example, Xie and Siegmund (2013) proposed to use the log-likelihood with a reference value for the changed proportion and Wu (2019) used the sum of S-R processes as the detection procedure when the proportion is small. These techniques works well if the signal are uniform and stable on the changed streams. Without loss of generality, we shall assume $\Sigma =I_N$. So the regular MEWMA chart makes an alarm at
\begin{equation}
\tau_{MEW} =\inf\{ t>0: ~  Y_{t}^TY_t =\sum_{i=1}^N Y_{it}^2 >b^2 (\frac{\beta}{2-\beta}) \}.
\end{equation}

As we mainly focus on the performance of MEWMA procedure. There are three techniques that can be used. One is using a hard-threshold by trimming the lower values of $Y_{it}$ as its components directly measures the signal strength. The second is multiplying different weights as soft-threshold for $Y_{it}$ depending on its value. Third is the adaptive technique by directly changing the weight factor $\beta_i $ depending on the value of $Y_{i(t-1)}$ for each component. Here, we shall only consider the first two approaches as the adaptive procedure changes the nature of the process and makes the computation much more complicated.  For the first approach, by dropping the channels with $Y_{it}^2 \geq s $, say $s=0.5$, we make an alarm at
\begin{equation}
\tau_{TEW} =\inf\{t>0: \sum_{j=1}^N Y_{jt}^2I_{[|Y_{jt}|>0.5]} > c^2. \}
\end{equation}

For the second approach, we can adapt the method proposed in Siegmund, et al. (2011) and use a soft-threshold  weight factor $e^{Y_{jt}^2/2}/(q+ e^{Y_{jt}^2/2})$ where $q=(1-p)/p$ is the ratio of proportions of no changed with changed channels. For example, if we take $p=0.1$, then $q=9$ and an alarm is made at
\begin{equation}
\tau_{WEW} =\inf\{ t>0: \sum_{j=1}^N \frac{e^{Y_{jt}^2/2}}{9+ e^{Y_{jt}^2/2}} Y_{jt}^2> d^2 \}. 
\end{equation}

Similarly, we can also use the above two approaches for the MMA method. For the regular MMA, we make an alarm at
\[
\tau_{MMA} =\inf\{ t>0: \sum_{t-w+1}^t \bar{X}_{ij}^2 > h^2 \}. 
\]
While for the hard-threshold MMA, we make an alarm at
\[
\tau_{TMA} =\inf\{t>0: \sum_{t-w+1}^t \bar{X}_{ij}^2I_{[|\bar{X}_{ij}|>0.5]} >h^2 \}. 
\]

In Table 4, we conduct a simulation comparison between several charts: the regular MEWMA chart; the MEWMA chart with hard-threshold 0.5 ; the MEWMA chart with soft-threshold 9; the MMA with window width $w=20$; and MA chart with hard-threshold 0.5. All simulations are replicated 10,000 times. The alarm limits are determined by using the approximations. Listed are the simulated $ARL_0$, false alarm rate (FAR), and the SADDT with $\nu =100$. The major finding is that for changed proportion less than equal 10\% (K=1,2), the hard-threshold EWMA reduces SADDT more than the soft-threshold EWMA. For $K=3, 5$, although the soft-threshold EWMA performs better than the hard-threshold EWMA, however, both have no improvement compared with the regular EWMA chart. That means, the threshold method works well only for sparse signal. Also, we observe that EWMA chart performs better than the MA chart, particularly for small and large signal strength. 

Listed also include the chart based on the sum of $N$ S-R processes as considered in Wu (2019) that is local optimal. With same reference value $\delta =0.5$ for $\mu$, it makes an alarm at
\[
\tau_{SSR} =\inf\{ t>0: \sum_{i=1}^N R_{it} >B \}
\]
where $R_{it} =(1+R_{i(t-1)})\exp(\delta (X_{it} -\delta /2 ))$. For $N=20$ and $ARL_0 =1000$, the approximation from $ARL_0 \approx \frac{B}{N} \exp(\delta \rho_+)$ gives $B \approx 14945.83$ . 
For a comparison, the SADDT for the regular S-R procedure as considered in Table 3 is also given in the bracket. 

We see that for small signal strength, the Sum-S-R procedure performs better than the hard-threshold EWMA. When signal strength gets larger ($>1$), the EWMA chart performs better. 

\begin{table}
\begin{center}
Table 4.  Comparison of SADDT for $ARL_0 =1000$, $N=20$ and $\nu=100$ with K changed channels \\
\begin{tabular}{|c |c |c | c| c|c |c|c| }   \hline
$\mu$ & K & EWMA & EWMA(0.5) & MA(WTD)  & MA & MA(0.5) & Sum-S-R (S-R)\\
 &  & $\beta=0.05$ & $\beta=0.05$ & $\beta =005$ & $w=20$ &$w=20$ & $\delta =0.5$\\
 &   & $\frac{b^2\beta}{2-\beta} =$ &$c^2=$ &$d^2=$ & $h^2 =$&  $h^2=$ & B=\\
 &  &  1.07 & 0.39 & 0.115 & 2.1125 & 1.26 & 14945.8(747.3) \\ \hline
$ARL_0$ & &1020.5 & 1052.74& 1063.60 & 1048.96& 1018.50 &991.31(1002.0)\\ \hline
FAR &  &0.0704&  0.0707  & 0.0691   & 0.0764 &  0.0783  & 0.0526(0.0717)\\ \hline
0.25 & 1 &382.65 & 352.321& 380.31  &578.39&  550.85 &182.71(473.22)\\
     & 2 &204.39 &207.07  & 203.98 & 354.89& 353.63 & 113.02(249.94) \\  
     & 3 &131.89& 148.28 &131.60   & 240.50&246.04 & 88.38(148.64)\\
     & 5 & 74.02 & 99.30  & 75.57 & 127.91 &137.52 & 65.29(68.82)\\ \hline
0.5 & 1&  93.65 & 66.14 & 86.11 & 172.78 & 155.59 & 48.18(251.82) \\
     & 2& 45.94 &41.83  & 44.55 & 68.09  & 60.07 & 35.11(96.00))\\  
     & 3 & 31.80 &32.92  &31.25   & 38.68 & 36.33 & 29.83(49.54)\\
     & 5 & 21.15 &25.20  &21.18   & 21.74 & 22.28 &24.11(23.32)\\ \hline
0.75& 1 &41.22 & 29.65& 37.24 &56.85  & 50.45& 26.55(149.05) \\
     & 2& 22.81 & 20.94 & 22.02 & 24.09  & 22.17 & 20.69(49.70)) \\  
     & 3 & 17.01 & 17.50 &16.73   &17.14  & 16.67 & 17.91(27.32) \\
     & 5 & 12.15 & 14.40 & 12.12  & 13.01 & 13.25 &15.06(13.70) \\ \hline
1.0 &  1 & 25.09 & 18.83& 22.88 & 27.47 & 24.83 & 18.36(96.63)\\
     & 2& 15.19 & 14.01  &14.68  & 15.45  & 14.71 & 14.61(31.96)) \\  
     & 3 & 11.72&   11.96  &11.52  &12.72 &12.36 & 12.95(18.07) \\
     & 5 & 8.64 &  10.03   & 8.61 & 9.98 & 10.15 & 11.01(9.60) \\ \hline
1.25 &  1& 18.11& 13.86 & 16.61 & 18.13 & 17.16 & 14.07(67.40) \\
     &  2& 11.47 &10.62  & 11.00 & 12.37  & 11.99 &11.37(23.38)) \\  
     & 3 & 8.96& 9.19  &8.82  & 10.21 & 10.11 & 10.06(13.61) \\
     & 5 & 6.70&  7.81  &6.72  & 8.11 & 8.34 &8.69(7.55)  \\ \hline
1.5 &  1&  14.06& 11.04 & 12.91 & 14.49 & 14.00 & 11.46(50.02)\\ 
     & 2& 9.14 & 8.52 &8.89  & 10.45  & 10.09 &9.34(18.00) \\  
     & 3 & 7.33 & 7.49&7.18   & 8.68 & 8.53 & 8.38(10.85) \\
     & 5 & 5.54 & 6.41 & 5.54  & 6.85 & 7.03 & 7.24(6.19)\\ \hline
1.75&  1 & 11.56&  9.13  & 10.62    & 12.48    &  12.13 & 9.57(39.77)\\
     & 2& 7.69 & 7.12 &7.45  & 9.08  & 8.80  &7.97(14.89)\\  
     & 3 & 6.18 & 6.33 & 6.08  & 7.54 & 7.43 & 7.17 (9.16)\\
     & 5 & 4.74 & 5.44 & 4.73  & 5.97 & 6.12 & 6.21 (5.25)\\ \hline
2.00&  1&  9.86 & 7.81 & 9.14     & 11.01   & 10.76 &8.31(31.92)  \\
     &  2& 6.63 & 6.20 &6.45  & 7.99  & 7.76 &6.97(12.65)  \\  
     & 3 & 5.37& 5.49 & 5.29  & 6.68 &6.58 & 6.24 (7.83)\\
     & 5 & 4.13 &4.76  & 4.15  & 5.27 & 5.41 & 5.44(4.61) \\\hline
\end{tabular}
\end{center}
\end{table}

\subsection{By adaptive estimation} 

Since MEWMA and MMA  charts are not only direction free and also do not need the signal strength, the corresponding SADDT may be larger when some information is available on the post change mean. As they provide natural adaptive estimates for the post-change mean, they can be used in other charts with certain optimal properties such as CUSUM and S-R procedures. For example, if we treat $\hat{\mu}_t$ as certain adaptive estimate of $\mu$ up to time $t$, the adaptive CUSUM procedure can be defined to make an alarm at
\[
\tau_{ACS}=\inf\{t>0: W_t =\max (0, W_{t-1} + \hat{\mu}_{t-1}^T (X_t -\hat{\mu}_{t-1}/2) ) >c\},
\]
and the adaptive S-R procedure can be defined to make an alarm at 
\[
\tau_{ASR} =\inf\{ t>0: R_t =(1+R_t) \exp(\hat{\mu}_{t-1}^T (X_t -\hat{\mu}_{t-1}/2) ) > B \}
\]

In Table 5, we present some simulation comparisons between several charts for the same setup as in Table 4. Included are the regular MEWMA, the adaptive S-R procedure with $\hat{\mu}$ is taken as the EWMA estimate, the adaptive CUSUM procedure. In the last column, we also include the sum of adaptive individual S-R procedure for detecting possible sparse signals that makes an alarm at
\[
\tau_{ASSR} =\inf\{ t>0: \sum_{i=1}^N R_{it} >B\},
\]
where $R_{it} = (1+R_{i(t-1)}) \exp (\hat{\mu}_{i(t-1)} (X_{it} - \hat{\mu}_{i(t-1)}/2))$. 

First, we find that the adaptive CUSUM and S-R procedure do not perform as well as the regular MEWMA procedure. So they do not need to be recommended. Second, in the sparse signal case, unless the signal is really sparse where the adaptive-Sum-S-R procedure performs better, no additional techniques are necessary to reduce the SADDT. However, the false alarm probability for the adaptive-sum-S-R procedure is much smaller than other charts. In summary, when the signal is sparse, the MEWMA procedure with hard threshold should be recommended. 

\begin{table}
\begin{center}
Table 5.  Comparison of SADDT for $ARL_0 =1000$, $N=20$ and $\nu=100$ with Adaptive Estimate \\
{\small 
\begin{tabular}{|c |c |c | c|c| c| }   \hline
$\mu$ & K &EWMA &  A-S-R(EWMA) &ACUSUM(EWMA) &A-Sum-S-R(EWMA)\\ 
  &  &  & $B=545$ &$c=4.6$ & $B=15050$  \\ \hline
$ARL_0$ & & 1020.5 & 995.2&1021.3 & 997.46\\
FAR & &0.0704 & 0.0900 & 0.0829 & 0.0187\\ \hline 
0.25 & 1 & 382.65 & 425.76 &446.53 & 196.12 \\
    & 2 &204.39 &  244.81&259.88 &135.53\\
    & 3 &131.89 &  156.94&167.63 & 108.89\\
    & 5 &74.02 &  87.07 &91.26&83.64 \\
    & 10 &37.32 &  40.72 & 42.05 & 60.18\\ \hline
0.5 & 1 & 93.65 & 113.39 &118.42 & 57.34\\
    & 2 & 45.94 &  50.65& 53.02 &42.72 \\
    & 3 & 31.80 & 33.50 &35.26 & 36.71\\
    & 5 &21.15&  21.61 & 22.06 & 30.34 \\
    & 10 & 13.02 & 13.18 & 13.27& 24.13 \\ \hline 
0.75 & 1 &41.22 &  44.56 &46.76 & 31.41\\
     & 2 &22.81 & 23.51 & 24.21 & 24.55\\
     & 3 & 17.01 & 17.39 &17.38 & 21.65\\
     & 5 & 12.15 & 12.23& 12.25 & 18.48 \\
     & 10 &8.09 &  8.13 &8.10 & 15.10 \\\hline
1.0 & 1 & 25.09 & 25.77 &26.96 & 21.45 \\
    & 2 &15.19 &  15.20 & 15.51 &17.33 \\
    & 3 & 11.72 & 11.75 &11.72 & 15.45 \\
    & 5 &8.64 &  8.65 &8.62 &13.41\\  
    & 10 &5.91 &  6.05 &5.96& 11.05 \\\hline
1.25 & 1 & 18.11 & 18.29 &18.37& 16.30 \\
     & 2 & 11.47 &  11.38& 11.46 &13.47 \\
     & 3 &8.96 &  9.00 &8.96 & 12.06\\
     & 5 &6.70 &  6.84 &6.78 &10.54\\
     & 10 &4.67 &  4.91 &4.81 & 8.86\\ \hline
1.5 & 1 &14.06 &  14.13 &14.25 & 13.17 \\
    & 2  &9.14 &  9.22 & 9.21 &10.98\\
    & 3  & 7.33 & 7.39 &7.31 & 9.84 \\
    & 5 &5.54 &  5.72 &5.60 &8.75\\ 
    & 10 & 3.93 & 4.22 &4.09 & 7.39 \\\hline
1.75 & 1 &11.56 &  11.58& 11.62 & 11.13\\
     & 2 & 7.69 & 7.75 & 7.69 & 9.38\\
     & 3 &6.18 &  6.31 & 6.23 & 8.51\\
     & 5 &4.74 &  4.97 & 4.89 &7.54\\ 
     & 10 &3.40 &  3.70 &3.62 & 6.38\\ \hline
2.0 & 1 & 9.86 & 9.88 & 9.84 & 9.65 \\
    & 2 &6.63&  6.75& 6.69 & 8.18 \\
    & 3 & 5.37 & 5.57&5.46 & 7.46\\
    & 5 &4.13 &  4.41 & 4.31 &6.63\\ 
    & 10 &3.00 &  3.35 &3.26 & 5.65\\\hline 
\end{tabular}
}
\end{center}
\end{table}

\section{Discussion and Conclusions}
In this communication, we studied theoretically the performance of several direction-invariant multivariate charts in detecting the change of mean including the $ARL_0$ and $SADDT$. In particular, accurate approximations for $ARL_0$'s are obtained and can be used for design purpose. Optimal design for MEWMA chart is studied and compared with MA, CUSUM, GLRT, and S-R charts. When the signal is sparse, hard and soft threshold techniques are studied for reducing SADDT. Overall, we find that MEWMA should be recommended due to its convenience in design and computation. There are several topics we did not cover. One is the effect of serial correlation within each channel. Combining with cross-sectional dependence, one can consider more complicated models as discussion in 
Rosenberg (1973) and Pesaran (2006) where common correlated effects (CCE) is considered by assuming randomly changing means.  Recent reviews and developments on dispersion control can be seen in Acosta-Mejiai, et al.(1999), Adjadi et al. (2021), Chan and Zhang (2001), Ebadi et al. (2021), and Xie et al. (2021). Theoretical results for detecting the dispersion change in multivariate case will be communicated in a future work. 

\section{Appendix: Proof of Theorem 5}

We only give the main steps. 
First we can write 
\[
\tau_{GLR} =\inf\{t>0: \max_{1\leq w \leq W} \sum_{j=1}^N Z_{jt}^{w2} > h^2 w =b^2\},
\]
where $Z_{jt}^w = \frac{1}{\sqrt{w}} \sum_{t=t-w+1}^t X_{ji} $ which follows $N(0,1)$ under $P_0^*(.)$. Note that 
\[
\psi(\theta) = \ln E_0^*[\exp(\theta Z_{jt}^{w2} )] =\frac{1}{2}\ln (1-2\theta).
\]
Define a changed measure $P_t^w(.) $ on the time window $[t-w+1, t] $ such that
$X_{t-w-1}, ..., X_t$ are i.i.d. $N(0, 1/\sqrt{1-2\theta})$ for $\theta <1/2$ with changed variance instead of mean. This implies that $Z_{it}^{w2}$ follows $\chi_1^2$ and the likelihood ratio for $X_1,,..., X_L$ has the form
\[
\frac{dP_t^w}{dP_{\infty}} (Z_{jt}^{w2}) = \exp(\theta  Z_{jt}^{w2} +\frac{1}{2} \ln(1-2\theta )),
\]
for $j=1,2.,,.N$. 
Let $\theta $ satisfy $N \frac{d}{d\theta} [-\ln(1-2\theta)/2] =b^2$, i.e. $1- 2\theta =N/b^2$ or $\theta = (1-N/b^2)/2 $. Thus $\theta $ is free of $w$.  By denoting
\[
l_t^w =\theta \sum_{j=1}^N Z_{jt}^{w2}+\frac{N}{2} \ln(1-2\theta );
\]
\[
\tilde{l}_t^w =\theta ( \sum_{j=1}^N Z_{jt}^{w2} -N/(1-2\theta)) ,
\]
and noting $\sigma_N^2 = Var_t^w(\tilde{l}_N ) = 2N\theta^2/(1-2\theta)^2$, for $L \rightarrow \infty $ and $W/L \rightarrow 0$, we can approximate 
\begin{eqnarray*}
P_{\infty}(\tau_{GLR} \leq L ) &=& P_{\infty} (\max_{1\leq t \leq L} \max_{1\leq w \leq W} \sum_{j=1}^N Z_{jt}^{w2} >b^2) \\
&=&\sum_{t=1}^L \sum_{w=1}^W E_{\infty}[ \frac{exp(l_t^w)}{\sum_{s=1}^L\sum_{1\leq u\leq W} exp(l_s^u)}; \max_{1\leq t\leq L} \max_{1\leq w \leq W}\sum_{j=1}^N Z_{jt}^{w2} >b^2] \\
&=& \sum_{t=1}^L\sum_{w=1}^W  E_t^w[ \frac{1}{\sum_{s=1}^L\sum_{1\leq u\leq W} exp(l_s^u)}; \max_{1\leq t\leq L} \max_{1\leq w \leq W}\sum_{j=1}^N Z_{jt}^{w2} >b^2 ] \\
&=& e^{-(\tilde{l}_t^w -l_t^w)} \sum_{t=1}^L \sum_{w=1}^W E_t^w[ \frac{M_t^w}{S_t^w} e^{-(\tilde{l}_t^w +\max_{1\leq s \leq L}\max_{1\leq u\leq W} (l_s^u-l_t^w))}; \tilde{l}_t^w + \max_{1\leq s \leq L}\max_{1\leq u\leq W} (l_s^u-l_t^w))]\geq 0]\\ 
&\approx& Le^{-N(\theta/(1-2\theta) +\ln(1-2\theta)/2)} \frac{1}{\sqrt{2\pi}}\frac{1-2\theta}{\sqrt{2N} \theta }
\sum_{w=1}^W \lim_{L\rightarrow \infty} E_t^w[ \frac{M_t^w}{S_t^w}]
\end{eqnarray*}
where $M_t^w =\max_{1<s\leq L} \max_{1\leq u\leq W}\exp (l_s^u-l_t^w) $ and $S_t^w = \sum _{w<s\leq L} \sum_{1\leq u\leq W} \exp (l_s^u-l_t^w) $. 
By using Sterling's formula for large $z$
\[
\Gamma(z) =\sqrt{\frac{2\pi}{z}} (\frac{z}{e})^z (1+O(\frac{1}{z})),
\]
we can show that 
\[
 e^{-N(\theta/(1-2\theta) +\ln(1-2\theta)/2} \frac{1}{\sqrt{2\pi}}\frac{1-2\theta}{\sqrt{2N} \theta }
 \approx \frac{(b^2/2)^{N/2-1} e^{-b^2/2}}{\Gamma(N/2) (1-N/b^2)}.
 \]
 
 Finally, we evaluate the factor $E_t^w[ \frac{M_t^w}{S_t^w}]$ by showing that $l_s^u-l_t^w$ can be approximated locally as two two-sided random walks. For technical convenience, we assume $N\rightarrow \infty $ and $N/b^2 \rightarrow 0$.   For $t-w \leq s \leq t$ and $u$ not too far way from $w$, we can write
 \begin{eqnarray*}
 Z_{js}^{u}-Z_{jt}^{w} & = & \frac{1}{\sqrt{u}} (\sum_{s-u+1}^s X_{ji} -\sum_{t-w+1}^t X_{ji}) +Z_{jt}^w (\sqrt{w/u}-1)  \\
 &\approx & \frac{1}{w} (\sum_{s-u+1}^{t-w}X_{ji} -\sum_{s+1}^t X_{ji} ) + Z_{jt}^w \frac{w-u}{2u}.
 \end{eqnarray*}
By ignoring the lower order terms, Taylor series expansion gives
\begin{eqnarray*}
Z_{js}^{u2}-Z_{jt}^{w2} & = & 2 Z_{jt}^{w} (Z_{js}^{u}-Z_{jt}^{w}) +(Z_{js}^{u}-Z_{jt}^{w})^2 \\
&\approx & \frac{2 Z_{jt}^{w}}{\sqrt{u}}(\sum_{s-u+1}^{t-w}X_{ji} -\sum_{s+1}^t X_{ji} ) +Z_{jt}^{w2} \frac{w-u}{u} +\frac{1}{u} (\sum_{s-w+1}^{t-w}X_{ji}^2 + \sum_{s+1}^t X_{ji}^2) \\
&=& \frac{2}{\sqrt{u}} \sum_{s-u+1}^{t-w} Z_{jt}^w X_{ji} +\frac{1}{u} \sum_{s-u+1}^{t-w} (X_{ji}^2-Z_{jt}^{w2}) - \frac{2}{\sqrt{u}} \sum_{s+1}^tZ_{jt}^w X_{ji} + \frac{1}{u} \sum_{s+1}^t 
(X_{ji}^2 +Z_{jt}^{w2}).
\end{eqnarray*}
Thus, we can write
\begin{eqnarray}
l_s^u -l_t^w &=& \theta \sum_{j=1}^N( Z_{js}^{u2}-Z_{jt}^{w2})  \nonumber\\
&\approx & \theta \frac{2\sqrt{N}}{\sqrt{u}} \sum_{s-u+1}^{t-w}  \frac{1}{\sqrt{N}} \sum_{j=1}^N Z_{jt}^w X_{ji}
+ \theta \frac{N}{u}\sum_{s-u+1}^{t-w} \frac{1}{N} \sum_{j=1}^N(X_{ji}^2-Z_{jt}^{w2})\\
& & - \theta \frac{2\sqrt{N}}{\sqrt{u}} \sum_{s+1}^{t} \frac{1}{\sqrt{N}}\sum_{j=1}^NZ_{jt}^w X_{ji} + \theta \frac{N}{u}
\sum_{s+1}^{t}\frac{1}{\sqrt{N}}\sum_{j=1}^N (X_{ji}^2 +Z_{jt}^{w2}) 
\end{eqnarray}

(7.18) and (7.19)  on the right side of the above equation define a two dimensional two-sided random walk. Follows the argument given in Siegmund, et al. (2011),  the random walk in (7.18) has drift coming from the second term as
\[
\theta \frac{N}{u} (1-\frac{b^2}{N}) = -\frac{b^2}{2u} (1- \frac{N}{b^2})^2 \approx -\frac{b^2}{2u},
\]
and variance coming from the first term as
\[
\theta^2 \frac{4N}{u} \frac{b^2}{N} = \frac{b^2}{u} (1- \frac{N}{b^2})^2 \approx \frac{b^2}{u},
\]
as $Z_{jt}^w$ and $X_{ji} $ are independent for $i\leq t-w$. Similarly, by using the conditional argument, we can show that the random walk defined by (7.19) has the same drift and variance. Same results are true for $t<s<t+w$. 
 From Siegmund and Yakir (2000), we can approximate
\[
 E_t^*[ \frac{M_t}{S_t}] \approx \frac{b^4}{4w^2} \upsilon^2 ( b/\sqrt{w}) \approx \frac{b^4}{4w^2} e^{-2 \rho_+ b/\sqrt{w}}.
 \]
 Thus, by approximating the summation by an integration, we have
 \begin{eqnarray}
 P_{\infty}(\tau_{GLR} \leq L ) & \approx & L \frac{(b^2/2)^{N/2-1} e^{-b^2/2}}{\Gamma(N/2)}
\int_0^W \frac{b^2}{4w} \upsilon^2 ( b/\sqrt{w}) dw \\
&\approx & Lb^2 \frac{(b^2/2)^{N/2-1} e^{-b^2/2}}{\Gamma(N/2)} \int_{b/\sqrt{W}}^{\infty} \frac{u \upsilon^2 (u) }{2} du.
\end{eqnarray}
By using the exponential approximation, we get the expected result for $ARL_0$. 
From this, we have
\[
b = \sqrt{2\ln(ARL_0)} (1+o(1)).
\]
At the alarm time $t >\nu$, with signal strength $||\mu||$, we have
\begin{eqnarray*}
& & \max w(\bar{X}_{t;w}^T\Sigma^{-1}\bar{X}_{t;w})\\
& & =\max  w((Z_{t,w}/\sqrt{w} +\min(\frac{t-\nu}{w},1)\mu)^T\Sigma^{-1} (Z_{t,w}/\sqrt{w} +\min(\frac{t-\nu}{w},1)\mu))\\
& & = \max w((\min(\frac{t-\nu}{w},1))^2 \mu^T\Sigma^{-1}\mu +2 \min(\frac{t-\nu}{w},1) 
\frac{\mu^T\Sigma^{-1}Z_{t;w}}{\sqrt{w}} + \frac{Z_{t;w}^T\Sigma^{-1}Z_{t;w}}{w}) =b^2. 
\end{eqnarray*}
Thus, as $b\rightarrow \infty$, $w\rightarrow \infty$. So the second and third term inside the square root are at lower order. Therefore, at the first order, we have
\[
\max \sqrt{w} \min (\frac{t-\nu}{w},1) ||\mu|| =b.
\]
Obviously, the maximum value $w$ satisfies $t-\nu = w$, i.e. $\sqrt{t-\nu} =b/||\mu|| +o_p(1).$. 
This gives 
\[
SADDT =b^2 /||\mu|^2 (1+o(1)) =2 \ln (ARL_0)/ ||\mu||^2 (1+o(1)). \qed
\]

\section*{Acknowledgement} Part of the work was done when the first author was visiting the Department of Statistics at Stanford University. 

\section*{References}
\begin{description}
\item[{\rm Adjadi, J. O. , Wang, Z., and Zwetsloot, I. M. (2021).}] 
A review of dispersion control charts for multivariate individual observations. {\em Quality Engineering} {\bf 33(1)}, 60-75.
\item[{\rm Acosta-Mejiai, J. J., Pignatiello, Jr., J. J., and Rao, B. V. (1999).}] A comparison of control charting procedures for monitoring
process dispersion. {\em IIE Transactions} {\bf 31}, 569-579
\item[{\rm  Bai, J. (2010).}] Common breaks in means and variances for panel data. {\em Journal of Econometrics} {\bf 157
(1)}, 78–92.
\item[{\rm Bardwell, L., Fearnhead, P., Eckley, I. Smith, S., and  Spott, M. (2019).}] 
Most recent changepoint detection in Panel data. {\em Technometrics} {\bf 61(1)}, 88-98
\item[{\rm Basseville, M. and Nikiforov, I. V. (1993).}] {\em Detection of Abrupt Changes: Theory and Applications}.  Prentice-Hall, Inc. 
\item [{\rm Chan, H. P. (2016).}]  Optimal sequential detection in multi-stream data. {\em Annals of Statistics} \textbf{45(6)}, 2736-2763.
\item[{\rm Chan, L.K. and Zhang, J. (2001).}] Cumulative sum control charts for the covariance matrix. {\em Statistica Sinica} {\bf 11}, 767-790. 
\item[{\rm Chen, J., Zhang W., and Poor, H.V. (2020).}] A false discovery rate oriented approach to
parallel sequential change detection problems.  {\em IEEE Transactions On Signal Processing} {\bf 68}, 1823-1836. 
\item[{\rm Chen, Y. and Li, X. (2019).}] Compound Sequential Change Point Detection in Multiple Data Streams. {\em Statistica Sinica}
\item[{\rm  Crosier, R.B. (1988).}] Multivariate generalizations of cumulative SUM quality control schemes.  {\em Technometrics} {\bf 30}, 291-303.
 \item[{\rm Ebadi, M., Chenouri, S., Lin, D.K.J., and  Steiner, S. H. (2021).}] Statistical monitoring of the covariance matrix in multivariate processes: a literature review. {\it Journal of Quality Technology} \\
 https://doi.org/10.1080/00224065.2021.1889419
 \item[{\rm Karlin, S. and Taylor, H.M. (1981).}] {\em A Second Course in Stochastic Processes}. Academic Press, New York.
\item[{\rm Kim, D. (2014).}]  Common breaks in time trends for large panel data with a factor structure. {\em Econometrics Journal } \textbf{17}, 301-337.
\item[{\rm Knoth, S. (2021).}] Steady-State average run length(s) - methodology, formulas and
numerics. {\em Sequential Analysis} {\bf 40(3)}, 405-426.
\item[{\rm Li, D., Qian, J., and Su, L. (2016).}] Panel Data Models With Interactive Fixed Effects and Multiple Structural Breaks. {\em Journal of American Statistical Association} \textbf{111}, 1804-1819. 
\item[{\rm Lowry, C.A., Woodall, W.H., Champ, C.W., and Rigdon, S.E. (1992).}] A multivariate exponentially weighted moving average control chart. {\em Technometrics}, {\bf 34(1)}, 46-53.
\item [{\rm  Mei, Y. (2010).}] Efficient scalable schemes for monitoring a large number of data
streams. {\em Biometrika} \textbf{97}, 419 -433.
\item[{\rm Ngai, H. M. and Zhang, J. (2001).}] Multivariate cumulative sum control charts based on projection pursuit. {\em Statistica Sinica} \textbf{11:} 747-766. 
\item[{\rm Nikiforov, I.V. (1999).}] Quadratic tests for detection of abrupt changes in multivariate signals. {\em IEEE Transactions on Signal Processing} {\bf 47(9)}, 2534-2538. 
\item[{\rm Ning, W. and Wu, Y. (2021).}] Common Change Point Estimation and Changed Panel Isolation after Sequential Detection in an Exponential Family. {\em Journal of Statistical Theory and Practice} {\bf 15(1)}\\
https://doi.org/10.1007/s42519-020-00134-3
 \item[{\rm Pesaran, M.H. (2006).}] Estimation and inference in large heterogeneous panels with a multifactor error structure.
 {\em Econometrika} \textbf{74(4)}, 967-1012.
 \item[{\rm  Pignatiello, J. J. and Runger, G.C. (1990).}]  Comparison of multivariate CUSUM charts {\em Journal of Quality Technology} \textbf{22}, 173-186.
 \item[{\rm Rosenberg, B. (1973).}] Random Coefficients Models: The analysis of a cross-sectional theory of time series by stochastically convergent parameter regression. {\em Journal of Economic and Social Measurements} \textbf{2(4)}, 399-428.
\item [{\rm Siegmund, D. (1985).}] {\it Sequential Analysis: Tests and Confidence Intervals}. Springer, New York.
\item[{\rm Siegmund, D. and Yakir, B. (2008).}] Detecting the emergence of a signal in a noisy image. {\em Statistics and its Interface} {\bf 1}, 3-12. 
\item[{\rm Siegmund, D., Yakir, B., and Zhang, N. R. (2011).}]
Detecting simultaneous variant intervals in aligned sequences. {\em Annals of Applied Statistics} {\bf 5(2A)}, 645-668. 
\item[{\rm Srivastava, M. S. and Wu, Y. (1993).}] Comparison of CUSUM, EWMA, and Shiryayev-Roberts procedures for detection of a shift in the mean. {\em Annals of Statistics} {\bf 21(2)}, 645-670. 
\item [{\rm Tartakovsky, A.G. and Veeravalli, V.V. (2008).}] Asymptotically optimal quickest
detection change detection in distributed sensor. {\em Sequential Analysis} \textbf{ 27}, 441-475.
\item[{\rm Woodall, W. H. and Ncube, M. M. (1985).}] Multivariate CUSUM Quality-Control Procedures. {\em Technometrics} \textbf{27(3)}, 285-292
\item[{\rm Wu, Y. (2019).}] A Combined SR-CUSUM Procedure for Detecting Common Changes in Panel Data. {\em Communication in Statistics: Theory and Methods} \textbf{48}(17), 4302--4319.
\item[{\rm Wu, Y. (2020).}] Estimation of common change point and isolation of changed panels after sequential detection. {\em Sequential Analysis} \textbf{39}(1), 52--64. 
\item[{\rm Wu, Y. and Wu, W. B. (2021).}] Sequential Common Change Detection, Isolation, and Estimation in Multiple Poisson Processes. To appear on {\em Sequential Analysis}.
\item[{\rm Wu, Y. and Wu, W. B. (2022).}] Sequential Detection of Common Transient Signals in High Dimensional Data Stream. {\em Naval Research Logistic Quarterly}. https://doi.org/10.1002/nav.22034
\item[{\rm Wu, Y. and Siegmund, D. (2022).}] Sequential Detection of Transient Signal in High Dimensional Data Stream. Submitted manuscript. 
\item [{\rm Xie, Y. and Siegmund, D. (2013).}] Sequential multi-sensor change-point detection.
{\em Annals of Statistics} \textbf{41}, 670-692.
\item[{\rm Xie, L., Xie, Y., and Moutstakides, G.V. (2021).}] Sequential subspace change-point detection. {\em Sequential Analysis} {\bf 39(3)},  307–335.
\end{description}

\end{document}